\documentclass[12pt]{article}
\usepackage{amsfonts}
\usepackage{amssymb}
\usepackage{amsmath}

\def\sqr#1#2{{\vcenter{\vbox{\hrule height.#2pt
              \hbox{\vrule width.#2pt height#1pt \kern#1pt \vrule width.#2pt}
              \hrule height.#2pt}}}}
\def\signed #1{{\unskip\nobreak\hfil\penalty50
              \hskip2em\hbox{}\nobreak\hfil#1
              \parfillskip=0pt \finalhyphendemerits=0 \par}}
\def\endpf{\signed {$\sqr69$}}

\def\dbR{{\mathop{\rm l\negthinspace R}}}
\def\dbC{{\mathop{\rm l\negthinspace\negthinspace\negthinspace C}}}

\def\3n{\negthinspace \negthinspace \negthinspace }
\def\2n{\negthinspace \negthinspace }
\def\1n{\negthinspace }

\def\dbC{{\mathop{\rm l\negthinspace\negthinspace\negthinspace C}}}

\def\ds{\displaystyle}

\def\dbN{{\mathop{\rm l\negthinspace N}}}

\def\dbR{{\mathop{\rm l\negthinspace R}}}
\def\={\buildrel \triangle \over =}

\def\a{\alpha}
\def\b{\beta}
\def\g{\gamma}

\def\e{\varepsilon}

\def\l{\lambda}

 \def\n{\nabla}

\def\t{\times}

\def\th{\theta}
\def\o{\omega}

\def\i{\infty}
\def\ns{\noalign{\ss} }

\def\G{\Gamma}
\def\D{\Delta}

\def\Si{\Sigma}

\def\O{\Omega}

\def\cA{{\cal A}}

\def\cL{{\cal L}}

\def\cO{{\cal O}}
\def\cP{{\cal P}}

\def\cR{{\cal R}}

\def\cl{{\cal l}}
\def\no{\noindent}

\def\ss{\smallskip}
\def\ms{\medskip}

\def\q{\quad}
\def\qq{\qquad}
\def\hb{\hbox}

\def\max{\mathop{\rm max}}

\def\pa{\partial}

\def\cd{\cdot}

\def\div{\hbox{\rm div$\,$}}

\def\cl{\overline}

\def\Re{{\mathop{\rm Re}\,}}
\def\Im{{\mathop{\rm Im}\,}}

\def\|{\Big |}
\def\({\Big (}
\def\){\Big )}
\def\[{\Big[}
\def\]{\Big]}
\def\be{\begin{equation}}
\def\bel{\begin{equation}\label}
\def\ee{\end{equation}}
\def\bt{\begin{theorem}}
\def\bcd{\begin{condition}}
\def\ecd{\end{condition}}
\def\et{\end{theorem}}
\def\bc{\begin{corollary}}
\def\ec{\end{corollary}}
\def\bde{\begin{definition}}
\def\ede{\end{definition}}
\def\bl{\begin{lemma}}
\def\el{\end{lemma}}
\def\bp{\begin{proposition}}
\def\ep{\end{proposition}}
\def\br{\begin{remark}}
\def\er{\end{remark}}
\def\ba{\begin{array}}
\def\ea{\end{array}}
\def\ed{\end{document}}
\def\ns{\noalign{\ms}}
\def\ds{\displaystyle}

\def\square#1{\vbox{\hrule\hbox{\vrule height#1%
     \kern#1\vrule}\hrule}}
\def\rectangle#1#2{\vbox{\hrule\hbox{\vrule height#1%
     \kern#2\vrule}\hrule}}

\font\tenbb=msbm10 \font\sevenbb=msbm7 \font\fivebb=msbm5

\newfam\bbfam
\scriptscriptfont\bbfam=\fivebb \textfont\bbfam=\tenbb
\scriptfont\bbfam=\sevenbb

\def\ov{{\overline v}}

\newtheorem{lemma}{Lemma}[section]
\newtheorem{remark}{Remark}[section]

\newtheorem{theorem}{Theorem}[section]
\newtheorem{corollary}{Corollary}[section]

\newtheorem{definition}{Definition}[section]
\newtheorem{proposition}{Proposition}[section]
\newtheorem{condition}{Condition}[section]

\makeatletter
   
   \@addtoreset{equation}{section}
\makeatother

\begin{document}
\title{\bf Stabilization of  the  weakly coupled
wave-plate system with  one internal
damping\thanks{This work is partially supported
by the  NSFC under grants 11471231, 11231007 and
11322110, by the Program for Changjiang Scholars
and Innovative Research Team in University under
grant IRT\_16R53, and the Fundamental Research
Funds for the Central Universities in China
under grant 2015SCU04A02.}}

\author{Xiaoyu Fu\thanks{School of Mathematics,  Sichuan University,
Chengdu 610064,  China. E-mail address:
rj\_xy@163.com.}\q and\q Qi L\"u\thanks{School
of Mathematics, Sichuan University, Chengdu
610064, China. E-mail address: lu@scu.edu.cn.}}

\date{}

\maketitle

\begin{abstract}

\no This paper is addressed to a stabilization
problem of a system coupled by a wave and a
Euler-Bernoulli plate equation. Only one
equation is supposed to be damped with a damping
function $d(\cd)$. Under some assumption about
the damping and the coupling terms, it is shown
that sufficiently smooth solutions of the system
decay logarithmically at infinity without any
geometric conditions on the effective damping
domain. The proofs of these decay results rely
on the interpolation inequalities for the
coupled elliptic-parabolic systems and make use
of the estimate of the resolvent operator for
the coupled  system. The main tools to derive
the desired interpolation inequalities are
global Carleman estimates.

\end{abstract}

\ms

\no{\bf Key Words}.  Logarithmic stability,
coupled wave-plate equations, interpolation
inequality, resolvent operator estimate.

\section{Introduction}

Let  $\Omega$ be a bounded domain in $\dbR^n$
($n\in \dbN$) with the $C^4$ boundary $\Gamma$.
Consider the following weakly coupled wave-plate system:
 \bel{0a1}\left\{\ba{ll}\ds
 y_{tt}-\D y+c(x) z+\a d(x)y_t=0 &\hb{ in } \dbR^+\t \O,\\
 \ns\ds
 z_{tt}+\D^2 z+c(x) y+(1-\a)d(x)z_t=0 &\hb{ in } \dbR^+\t \O,\\
 \ns\ds
y=z=\D z =0&\hb{ on } \dbR^+\t \G,\\
\ns\ds (y(0), y_t(0))=(y^0,y^1),\ (z(0), z_t(0))=(z^0,z^1)&\hb{ in }
\O.
 \ea\right.\ee
Here $c(\cd)\in L^\infty(\O)$ denotes the
coupling function, $d(\cd)\in L^\infty(\O)$
denotes the damping function, $\a=0$ or $\a=1$.
Both $c(\cd)$ and $d(\cd)$ are nonnegative
functions.

\ms

When $d=0$, the system \eqref{0a1} is a
classical model to describe the propagation of
waves in elastic solids
(e.g.\cite{Achenbach,Graff}). In recent years,
the model in which $d\neq 0$  attracts lots of
attentions due to the study of the plastic
composite materials. These materials are widely
used in industry, such as aerocraft, ships,
submarines and automobiles. Light weight is one
of the main advantage of them. However, light
weight could lead the structural elements of the
composite to subject to unnecessary vibrations.
Then it is important to add some damping on the
system to attenuate unnecessary vibrations in
the design of composite dynamic structures(e.g.
\cite{Adams,Yim}). In \eqref{0a1}, if $\a=1$,
then the term ``$d(x)y_t$" is an damping acted
on the wave which is used to stabilize the wave
directly and the plate indirectly. On the other
hand, if $\a=0$, then ``$d(x)z_t$" is an damping
acted on the plate which is used to stabilize
the plate directly and the wave indirectly.

The stabilization of both wave and plate
equations are studied extensively in the
literature (see
\cite{BLR,BH,Fu1,L,LR,Ortega,P1,Zua} and the
rich references therein). Generally speaking,
there are three types of decays for the energy
of damped systems, that is, exponential decay,
polynomial decay and logarithmic decay.  To
obtain the first two decays, one needs
restrictive conditions on the support of the
damping. For example, to get the exponential
decay of the energy of wave equations, one needs
the set
$$
\o\=\{x\in\O:\, d(x)\geq d_0\mbox{ for some
$d_0>0$}\}
$$
fulfills the geometric control condition (e.g.
\cite{BLR}). To get the polynomial decay, one
needs $\O$ and $\o$ fulfills some special
geometric condition (see \cite{BH,BZ,P1} for
example). Similar things happen for plate
equations (e.g.\cite{ALM,BZ1,Haraux}).

When $\O$ and $\o$ does not fulfill any special
condition, people find that the energies of both
wave and plate equations satisfy logarithmic
decay(e.g. \cite{L,LR}). In this paper, we show
such decay also holds for \eqref{0a1}.

There exist many results on the stabilization of
coupled systems in the literature (see
\cite{AB99, B1,ACK,AN,BL,RZZ,ZZ2} and the rich
references therein). Particularly, in
\cite{ACK}, the author obtained the polynomial
decay of the system \eqref{0a1} with a constant
coupling parameter and a damping which is
effective on the whole boundary. As far as we
know, there is no reference addressing the
asymptotic behavior of the system (\ref{0a1})
(with only one damping but without geometric
assumptions on the effective damping domain).

In this paper, we will show the logarithmic
decay property for solutions of the system
(\ref{0a1}).  According to a well known result
of Burq (see \cite[Theorem 3]{B}),   to obtain
the logarithmic decay rate, it suffices to show
some high-frequency estimates with exponential
loss on the resolvent. Thus, the main difficulty
is the estimate of the resolvent operators,
which will be solved by using some Carleman
inequalities. To this end, we borrow some idea
in \cite{Fu3}. However, there are some new
essential difficulties. Indeed, to get the
energy decay for a system coupled by two wave
equations, one should establish an interpolation
inequality for a system coupled by two elliptic
equations. In our case, we have to get two
interpolation inequalities for a system coupled
by one elliptic and two parabolic equations. One
cannot simply mimic the techniques in \cite{Fu3}
to obtain our result. Please see Section
\ref{ss5} for more details.

The rest of this paper is organized as follows:
In Section \ref{ops}, we give the main results
in this paper. Section \ref{sec-Carleman}  is
devoted to establish Carleman estimates for some
second order partial differential operators.
Section \ref{ss5} is addressed to proving  some
interpolation inequalities by those Carleman
estimates. At last, in Section \ref{ss6}, we
prove our main results.

\section{Statement of the main results}\label{ops}

Throughout this paper, we always assume that
$c(\cd)$ and $d(\cd)$ are bounded real valued
nonnegative functions satisfying
 \bel{cond1}
c(x)\ge c_0>0  \ \hbox{ in } \o_c \ee and
\bel{cond1.1}d(x)\ge d_0>0 \ \hbox{ in }\o_{d},
\ee where $\o_c$ and $ \o_d$ are any fixed
non-empty open subsets of $\O$, and $c_0$ and
$d_0$ are given constants.

\begin{remark}
The condition \eqref{cond1} means that the two
equations are at least coupled effectively on
the non-empty open set $\o_c$. If it does not
hold, then the wave and the plate may not impact
each other efficaciously. Then one cannot
stabilize the system by only put damping on one
equation.
\end{remark}
\begin{remark}
The condition \eqref{cond1.1} means that the
damping are at least worked effectively on the
non-empty open set $\o_d$. If it does not hold,
then the damping may not be strong enough to
stabilize the system. An extreme example is that
$d=0$.
\end{remark}

In what follows, we will use $\ds
C=C(\O,\o_c,\o_{d})$ to denote generic positive
constants which may vary from line to line.


Let
$$
 H\=H_0^1(\O)\t L^2(\O)\t (H^2(\O)\cap H_0^1(\O))\t L^2(\O).
 $$
Write an element $U\in H$ as
  $
  U=(y,u,z,v)$, where $ y\in H_0^1(\O)$ and $ u, v\in
  L^2(\O),\ z\in H^2(\O)\cap H_0^1(\O)$. Define an unbounded operator $\cA: \; D(\cA)\subset H\to
H$ by
 \bel{0a2}\left\{\ba{ll}
 D(\cA)\=\left\{U\in H: \  \cA U\in H, \ y\big|_{\G}=z\big|_{\G}=\D z\big|_{\G}=0\right\}, \\
 \ns\ds \cA U=\(u, \D y- c(x) z-\a d(x)u,v, -\D^2 z-c(x)y-(1-\a)d(x)v\).\ea\right.\ee
It is easy to show that $\cA$ generates a
$C_0$-semigroup $\{e^{t\cA}\}_{t\geq 0}$ on $H$.
Therefore, system (\ref{0a1}) is well-posed. The
energy of a solution $(y,z)$ to the system
(\ref{0a1}) at time $t$ is given by:
 \bel{energy}\ba{ll}\ds
E(t)&\ds=\frac{1}{2}\int_\O\big(|\n
y(t)|^2+|y_t(t)|^2\big)dx\\
\ns&\ds\q+\frac{1}{2}\int_\O\big(|\D
z(t)|^2+|z_t(t)|^2\big)dx+\int_\O c(x)\Re\(
y(t)\overline z(t)\)dx.
 \ea\ee

When $d=0$, it is easy to see that $E(\cd)$ is
conservative. When $d\neq 0$, our main results
are stated as follows:
 \bt\label{0t1}
Let
 $c(\cd)$ and $d(\cd)$  satisfy (\ref{cond1}) and (\ref{cond1.1}). Suppose that
$\o_c\cap\o_{d}\neq\emptyset$. Then, solutions
$e^{t\cA}(y^0,y^1,z^0,z^1)\equiv
(y,y_t,z,z_t)\in C(\dbR;\; D(\cA))\cap
C^1(\dbR;\;H)$ to system (\ref{0a1}) satisfy
  \bel{0a3}\ba{ll}\ds
 ||e^{t\cA}(y^0,y^1,z^0,z^1)||_{H}\le \frac{C}{\log (2+t)}||(y^0,y^1,z^0,z^1)||_{D(\cA)},\\
 \ns\ds\qq\qq\qq\qq\qq\qq\forall\;(y^0,y^1,z^0,z^1)\in D(\cA),\ \forall\; t>0.
  \ea\ee
 \et
\begin{remark}
From \cite{L}, we know the decay rate in Theorem
\ref{0t1} is sharp.
\end{remark}

Following \cite[Theorem 7.1]{D} (see also
\cite[Theorem 3]{B}), Theorem \ref{0t1} is a
consequence of the following resolvent estimate
for the operator $\cA$:

 \bt\label{0t2} Under the assumptions of Theorem \ref{0t1}, there
exists a constant $C>0$ such that if
 $$
 \Re\l\in\left[-{e^{-C|\Im\l|}/C},0\right],
 $$
then it holds that
  $$
 ||(\cA-\l I)^{-1}||_{\cL(H)}\le Ce^{C|\Im\l|}, \q\hb{ for } |\l|>1.
  $$
 \et

 \br
In this paper, we assume that
$\o_c\cap\o_{d}\neq\emptyset$. It would be quite
interesting to consider the case that
$\o_c\cap\o_{d}=\emptyset$. As far as we know,
this is an unsolved problem.
 \er
%


\section{Carleman estimates for the  elliptic and
parabolic operators} \label{sec-Carleman}


In this section,  we establish Carleman
estimates for elliptic operator $\pa_{ss}+\D$,
and parabolic operators $\pa_t+\D$ and
$\pa_s-\D$, respectively.

To begin with, we assume that $\o_k$
($k=0, 1, 2, 3, 4$) are subdomains of $\O$ such that
$\o_0\subset\o_1\subset\o_2\subset\o_3\subset\o_c\cap\o_{d}\=\o_4$.
Recall that there exists a function $\hat\psi\in
C^2(\overline \O)$ such that (see \cite{FI} for
example)
 \bel{as}
 \hat\psi> 0  \hb{ in } \O,
 \q\hat\psi=0 \hb{ on }\pa\O,\q
|\n\hat\psi|>0  \hb{ in }\overline{\O\setminus\o_0}.
 \ee
With the aid of the function $\hat\psi$ defined above, we introduce
 weight functions as follows:
 \bel{as4}
 \th=e^{\ell}, \q \ell=\l\phi,\q \phi=e^{\mu\psi},\q \psi=\psi(s,x)\=\frac{\hat\psi(x)}{
 ||\hat\psi||_{L^\infty(\O)}}+b^2-s^2.
  \ee
Here $1<b\le2$ will be given later, $\l, \mu,
s\in \dbR$ are parameters and $x\in\cl{\O}$.

We first recall the  following known global Carleman estimate for elliptic
equations.
 \bl\label{car1} (\cite[Theorem 3.1]{Fu3})
Let  $p\in C^2((-b,b)\t\O;\;\dbC)$, and let
$\ell\in C^2((-b,b)\t\O)$ be given by
(\ref{as4}). Then, there is a constant $\mu_0>0$
such that for all $\mu\ge \mu_0$, one can find
two constants $C=C(\mu)>0$ and $\l_0=\l_0(\mu)$
so that for all $p\in H_0^1((-b,b)\t\O)$ and
$\ds p_{ss}+\D p=f_1$ (in $(-b,b)\t\O$, in the
sense of distribution) with $f_1\in
L^2((-b,b)\t\O)$, and for all $\l\ge\l_0$, it
holds that
 \bel{crf}\ba{ll}\ds
 \l\mu^2\int_{-b}^b\int_{\O}\th^2\phi(|\n p|^2+|p_s|^2+\l^2\mu^2\phi^2|p|^2)dxds\\
\ns\ds\le
C\Big[\int_{-b}^b\int_{\O}\th^2|f_1|^2dxds+\l\mu^2\int_{-b}^b\int_{\o_0}\th^2\phi(|\n
p|^2+|p_s|^2+\l^2\mu^2\phi^2|p|^2)dxds\Big].
 \ea\ee
 \el

 \ms

Further, choose some cut-off functions
$\eta_{j+1} \in C_0^\i(\o_{j+1})$ ($j=0, 1, 2,
3, 4$) such that $\eta_{j+1} (x)=1$ in $\o_{j}$.
Then, we have the following local weighted
energy estimate.

  \bl\label{le0m-01}
  Let $\g\in\dbR$ and $\ell\in C^2((-b,b)\t\O;\dbR)$ be given by (\ref{as4}). Then,
there is a constant $\l_0>0$ such that for all $\l\ge \l_0$, one
can find  a constant $C>0$   so that for
all $q\in H_0^1((-b,b)\t\O)$ and $\ds
\g q_{s}+\D q=f_2$ (in $(-b,b)\t\O$, in
the sense of distribution) with $f_2\in L^2((-b,b)\t\O)$, for any $\b\ge 2$, it holds that
 \bel{1224-x1}\ba{ll}\ds
 \int_{-b}^b\int_{\O}\th^2\eta_{j+1}^2\phi^k |\n q|^2dxds\\
 \ns\ds\le \frac{1}{(\l\mu)^\b}\int_{-b}^b\int_{\O}\th^2|f_2|^2dxds+C(\l\mu)^\b\int_{-b}^b\int_{\O}\th^2\eta_{j+1}^2\phi^{k+2}|q|^2
dxds
 \ea\ee
  where $j=0, 1, 2, 3, 4$ and $k\in\dbN$.
  \el

{\it Proof. } We choose a cut-off function
$\eta_{j+1} \in C_0^\i(\o_{j+1}; [0,1])$ such
that $\eta (x)=1$ in $\o_{j}$. A short
calculation shows that
  \bel{1202eq1}\ba{ll}\ds
\th^2\phi^k\eta_{j+1}^2\[\overline {q}(\g q_{s}+\D {q})+ {q}\overline{(\g q_{s}+\D {q})}\]\\
\ns\ds=(\g\th^2\phi^k\eta_{j+1}^2|q|^2)_s-(\g\th^2\phi^k\eta_{j+1}^2)_{s}|q|^2+\sum_{l=1}^n\[\th^2\phi^k\eta_{j+1}^2(\overline
{q}q_{x_l}+q{\overline q}_{x_l})\]_{x_l}\\
\ns\ds\qq-2\th^2\phi^k\eta_{j+1}^2|\n q|^2-\sum_{l=1}^n\(\th^2\phi^k\eta_{j+1}^2\)_{x_l}(\overline
{q}q_{x_l}+q{\overline q}_{x_l}).
  \ea
  \ee
Next,  recalling (\ref{as4}) for the definition
of $\ell$, it is easy to see that
 \bel{1202-as6}\left\{\ba{ll}\ds
\ell_s=\l\mu\phi\psi_s,\qq\ell_{x_j}=\l\mu\phi\psi_{x_j},\q\ell_{x_j s}=\l\mu^2\phi\psi_s\psi_{x_j},\\
\ns\ds
\ell_{ss}=\l\mu^2\phi\psi_s^2+\l\mu\phi\psi_{ss},\q\ell_{x_jx_k}=\l\mu^2\phi\psi_{x_j}\psi_{x_k}+\l\mu\phi\psi_{x_jx_k}.
 \ea\right.
 \ee

Integrating (\ref{1202eq1}) in $(-b,b)\t\O$,  noting that $q (-b)=q(b)=0$ in $\O$ and $q=0$ on the boundary, by (\ref{1202-as6}), we get the desire result immediately. \endpf

\ms

Further, we recall the following point-wise weighted inequality for  the parabolic operators.
 \bl\label{1218t1} (\cite[Theorem 2.1]{Fu01}) Let $\g\in \dbR$ and $z\in C^2(\dbR^{1+n}; \;\dbC)$. Set $\th=e^\ell,\ v=\th
q$ and $\ds\Psi=-2\D\ell$. Then
 \bel{1218a2}\ba{ll}\ds
 \frac{1}{2}\th^2|\g q_s+\D q|^2+M_s+\div V\\
\ns\ds\ge2
\sum_{j,k=1}^n\ell_{x_jx_k}(v_{x_k}\ov_{x_j}+\ov_{x_k}
v_{x_j})+2\D\ell|\n v|^2+B|v|^2,
 \ea\ee
 where
 \bel{1218f3}
 A\=|\n\ell|^2+\D\ell,
 \ee
 and
 \bel{1218a3}\left\{\ba{ll}\ds
M\=(\g^2 \ell_s-\g A) |v|^2+\g|\n v|^2,\\
\ns\ds V\=[V^1,\cdots,V^k,\cdots,V^n],\\
 \ns\ds V^k\=-\g (v_{x_k}\ov_s+\ov_{x_k}v_s)+2\sum_{j=1}^n\ell_{x_j}(v_{x_j}\ov_{x_k}+\ov_{x_j}v_{x_k})-2\ell_{x_k}|\n v|^2\\
\ns\ds\qq\q +2\D\ell(v_{x_k}\ov+\ov_{x_k}v)
+(2A\ell_{x_k}-2\D\ell_{x_k}-2\g\ell_{x_k}\ell_s)|v|^2,\\
\ns\ds
B\=\g^2\ell_{ss}-\g A_s-2\g(\n\ell\cdot\n\ell_s-\D\ell\ell_s)-2\D^2\ell+2(\n A\cdot\n \ell-A\D\ell).
 \ea\right.\ee \el

\ms

{\it Proof.}  Taking $\a=\g\in\dbR,\ \b=0$, $m=n$ and $(a^{jk})_{n\t n}= I_n$ (the unit matrix) in  \cite[Theorem 2.1]{Fu01} with $t$ replaced by $s$, in this case $\cP q=\g q_s+\D q$.
Noting that $\Psi=-2\D\ell$ and $\ds \th(\cP q\overline {I_1}+\overline{\cP q} I_1)\le 2|I_1|^2+\frac{1}{ 2}\th^2|\cP q|^2$, we immediately get the desired result.\endpf

 \br
Similar point-wise weighted identity was given in  \cite[Lemma 2.1]{lz} for the real valued parabolic equation.
 \er

 \ms

Based on Lemma \ref{le0m-01} and Lemma \ref{1218t1}, we have the following Carleman estimate  for the  heat equations.
 \bl\label{1202-car1}
Let $\g\in\dbR$ and $\ell\in C^2((-b,b)\t\O;\dbR)$ be given by (\ref{as4}). Then,
there is a constant $\mu_1>0$ such that for all $\mu\ge \mu_1$, one
can find two constants $C=C(\mu)>0$ and $\l_1=\l_1(\mu)$ so that for
all $q\in H_0^1((-b,b)\t\O)$ and $\ds
\g q_{s}+\D q=f_2$ (in $(-b,b)\t\O$, in
the sense of distribution) with $f_2\in L^2((-b,b)\t\O)$, and for all
$\l\ge\l_1$, it holds that
 \bel{1201-crf}\ba{ll}\ds
\int_{-b}^b\int_{\O}\th^2(\l\phi)^{-1}(|\D q|^2+|\g q_s|^2)dxds+
 \l\mu^2\int_{-b}^b\int_{\O}\th^2\phi(|\n q|^2+\l^2\mu^2\phi^2|q|^2)dxds\\
\ns\ds\le C\[\int_{-b}^b\int_{\O}\th^2|f_2|^2dxds+\l^3\mu^4\int_{-b}^b\int_{\o_1}\th^2\phi^3|q|^2dxds\].
 \ea\ee
 \el

\ms

{\it Proof. } We divide the proof into several steps.

\ms

{\it Step 1. } By (\ref{1202-as6}), a short calculation shows that
 \bel{1218x1}\ba{ll}\ds
 2\sum_{j,k=1}^n\ell_{x_jx_k}(v_{x_k}\ov_{x_j}+\ov_{x_k}
v_{x_j})+2\D\ell|\n v|^2\\
\ns\ds\ge 4\l\mu^2\phi |\n\psi\cdot\n v|^2+2\l\mu^2\phi |\n\psi|^2|\n v|^2-C\l\mu\phi |\n v|^2.
 \ea\ee
Next, by (\ref{as4}), (\ref{1202-as6}) and
(\ref{1218f3}), we obtain that
 \bel{1218x3}\left\{\ba{ll}\ds
 A=\l^2\mu^2\phi^2|\n\psi|^2+\l\mu^2\phi|\n\psi|^2+\l\mu\phi\D\psi\\
 \ns\ds \n (\phi^2|\n\psi|^2)\cdot\n\ell=2\l\mu^2\phi^3 |\n\psi|^4+\l\mu\phi^3\n(|\n\psi|^2)\cdot\n\psi.
 \ea\right.\ee
Recalling (\ref{1218a3}) for the definition of
$B$, it follows from (\ref{1218x3}) that
 \bel{1218x2}\ba{ll}\ds
 B\3n&\ds=\g^2\ell_{ss}-\g A_s-2\g(\n\ell\cdot\n\ell_s-\D\ell\ell_s)-2\D^2\ell+2(\n A\cdot\n \ell-A\D\ell)\\
 \ns&\ds\ge 2\l^3\mu^4\phi^3|\n\psi|^4-C\l^3\mu^3\phi^3-C\l\mu^4\phi.
 \ea\ee
Combining (\ref{1218x1}) and (\ref{1218x2}), we
find that
 \bel{0619-1}\ba{ll}\ds
 \hb{ RHS of (\ref{1218a2})}\3n&\ds\ge 2\l\mu^2\phi |\n\psi|^2 (|\n v|^2+\l^2\mu^2\phi^2|\n\psi|^2v^2)\\
 \ns&\ds\q -C\l\mu\phi(|\n v|^2+\l^2\mu^2\phi^2|v|^2+\mu^3 |v|^2).
 \ea\ee
Integrating inequality (\ref{1218a2}) on
$(-b,b)\t\O$, using integration by parts, noting
that $v_{x_j}=\frac{\pa v}{\pa\nu}\nu_j$ on
$\Si$ (which follows from $v|_\Si=0$),
$v(-b)=v(b)=0$ and by (\ref{as}) and
(\ref{0619-1}), we have that
 \bel{1202oo1}\ba{ll}\ds
\l\mu^2\int_{-b}^b\int_\O\phi |\n\psi|^2\(|\n
v|^2+\l^2\mu^2\phi^2|\n\psi|^2|v|^2\)dxds\\
\ns\ds\le C\[||\th f_2||^2_{L^2(Q)}+\l\mu\int_{-b}^b\int_\O\phi
\(|\n v|^2+\l^2\mu^2\phi^2|v|^2+\mu^3|v|^2\)dxds\],
 \ea\ee
 where we  use the following fact
  $$
  \int_{-b}^b\int_\O \div Vdxds=2\l\mu\int_{-b}^b\int_\G \phi \frac{\pa\psi}{\pa\nu}\|\frac{\pa v}{\pa\nu}\|^2dxds\le0.
  $$

On the other hand, by (\ref{as}) and (\ref{as4}), it is easy to see that
 $$
h\=|\n\psi|=\frac{1}{||\hat\psi||_{L^\infty(\O)}}|\n\hat\psi(x)|>0,\qq \hb{ in
}\cl{\O\setminus\o_0}.
 $$
Then
 \bel{1202oo4}\ba{ll}&\ds
\l\mu^2\int_{-b}^b\int_\O\phi h^2\(|\n
v|^2+\l^2\mu^2\phi^2h^2|v|^2\)dxds\\
\ns&\ds\ge c_1\l\mu^2\int_{-b}^b\int_{\O\setminus \O_{\o_0}}\phi
\(|\n
v|^2+\l^2\mu^2\phi^2|v|^2\)dxds\\
\ns&\ds\q-C\l\mu^2\int_{-b}^b\int_{{\o_0}}\phi \(|\n
v|^2+\l^2\mu^2\phi^2|v|^2\)dxds.
 \ea\ee
Combining (\ref{1202oo1})--(\ref{1202oo4}),
noting that $v=\th q$, we get that
 \bel{1224-0}\ba{ll}\ds
 \l\mu^2\int_{-b}^b\int_{\O}\th^2\phi\(|\n
q|^2+\l^2\mu^2\phi^2|q|^2\)dxds\\
\ns\ds\le C\[||\th f_2||^2_{L^2(Q)}+\l\mu\int_{-b}^b\int_\O\th^2\phi
\(|\n q|^2+\l^2\mu^2\phi^2|q|^2+\mu^3|q|^2\)dxds\]\\
\ns\ds\qq+C\l\mu^2\int_{-b}^b\int_{\o_0}\th^2\phi\(|\n
q|^2+\l^2\mu^2\phi^2|q|^2\)dxds.
 \ea\ee
Taking $\ds\mu_1\=4C+1>0$, for any
$\mu\ge\mu_1$, we have
 \bel{1224-1}\ba{ll}\ds
 \l\mu^2\int_{-b}^b\int_{\O}\th^2\phi\(|\n
q|^2+\l^2\mu^2\phi^2|q|^2\)dxds\\
\ns\ds\le C\(||\th f_2||^2_{L^2(Q)}+\l\mu^4\int_{-b}^b\int_\O\th^2\phi |q|^2dxds\)\\
\ns\ds\qq+C\l\mu^2\int_{-b}^b\int_{\o_0}\th^2\phi\(|\n
q|^2+\l^2\mu^2\phi^2|q|^2\)dxds.
 \ea\ee

Taking $\l_2\=2\sqrt{C}+1>0$, for any
$\l\ge\l_2$, we have
 \bel{611-0}\ba{ll}\ds
 \l\mu^2\int_{-b}^b\int_{\O}\th^2\phi(|\n q|^2+\l^2\mu^2\phi^2|q|^2)dxds\\
\ns\ds\le C\[\int_{-b}^b\int_{\O}\th^2|f_2|^2dxds+\l\mu^2\int_{-b}^b\int_{\o_0}\th^2\phi(|\n
q|^2+\l^2\mu^2\phi^2|q|^2)dxds\].
 \ea\ee

 \ms

 {\it Step 2. } Let us estimate ``$\ds\int_{-b}^b\int_{\O}\th^2(\l\phi)^{-1}(|\D q|^2+|\g q_s|^2)dxds$".

By  (\ref{1202-as6}) and (\ref{as}), a short calculation shows that
\begin{eqnarray}\label{1202-213}
\begin{array}{ll}
\displaystyle\int_{-b}^b\int_{\O}
\th^2 (\l\phi)^{-1}
\g (\overline q_s\cdot\D q+q_s\overline\D q)dxds\\
\ns\ds =-\g\displaystyle\int_{-b}^b\int_{\O}\displaystyle \n
\big(\th^2(\l\phi)^{-1}\big)\cdot(\n \overline q q_s+\n q \overline q_s)
dxds
+\g\displaystyle\int_{-b}^b\int_{\O}\displaystyle\big(\th^2(\l\phi)^{-1}\big)_s|\n q|^2 dxds\\
\ns\ds
\leq\displaystyle\frac{1}{2}\displaystyle\int_{-b}^b\int_{\O}
\th^2 (\l\phi)^{-1}|\g
q_s|^2dxdt+C\l\mu^2\int_{-b}^b\int_{\O}\th^2\phi |\n q|^2 dxds.
\end{array}
\end{eqnarray}

 Combining (\ref{1202-213}) and (\ref{611-0}), we end up with
\begin{eqnarray}\label{1202-214}
\begin{array}{ll}
\ds\int_{-b}^b\int_{\O}\th^2(\l\phi)^{-1}(|\D q|^2+|\g q_s|^2)dxds\\
\ns\ds \le C\[\int_{-b}^b\int_{\O}\th^2|f_2|^2dxds+\l\mu^2\int_{-b}^b\int_{\o_0}\th^2\phi(|\n
q|^2+\l^2\mu^2\phi^2|q|^2)dxds\].
\end{array}
\end{eqnarray}

Finally,
Combing (\ref{611-0}) and  (\ref{1202-214}),  by Lemma \ref{le0m-01}, taking $\l_1=\max\{\l_0,\ \l_2\}$, we get the desired estimate in
Lemma \ref{1202-car1} immediately. \endpf

\section{Two interpolation inequalities for the weakly coupled elliptic-parabolic system}\label{ss5}

In this section we
shall prove two interpolation inequalities for the following weakly coupled
elliptic-parabolic system:
  \bel{1202-s1}\left\{\ba{ll}\ds
p_{ss}+\D p-c(x)q+i\a d(x) p_s=p^0 & \mbox{ in } X, \\
\ns\ds w_{s}+\D w-c(x)p+i(1-\a)d(x) q_s=w^0 & \mbox{ in } X, \\
\ns\ds q_{s}-\D q=w & \mbox{ in } X, \\
\ns\ds p=q=w=0&\mbox{ on }\Si.
 \ea\right.\ee
  Here $X=(-2,2)\t\O,\ \Si=(-2,2)\t\G$ and $p^0,\ w^0\in
L^2(X)$. In what follows, we will use the notations
$Y\=(-1,1)\t\O,\q X^*\=(-2,2)\t(\o_c\cap\o_{d})$.
We have the following interpolation inequalities for system
(\ref{1202-s1}).
 \bt\label{1202-t1} Let $\a=1$. Under the assumptions in Theorem \ref{0t1},  there exists a constant $C>0$ such that, for any
$\e>0$, any
 solutions $(p,w, q)$ of system (\ref{1202-s1}) satisfies
  \bel{1202-0a2}\ba{ll}\ds ||p||_{H^1(Y)}+||w||_{L^2(Y)}+||q||_{L^2(Y)}\\
  \ns\ds\le
 Ce^{C/\e}\big(||p^0||_{L^2(X)}+||w^0||_{L^2(X)}+||p||_{L^2(X^*)}+||p_s||_{L^2(X^*)}\big)\\
 \ns\ds\q+Ce^{-2/\e}\big(||p||_{L^2(X)}+||p_s||_{L^2(X)}+||w||_{L^2(X)}+||q||_{L^2(X)}\big).
  \ea\ee
 \et
 \bt\label{1202-t2} Let $\a=0$. Under the assumptions in Theorem \ref{0t1},  there exists a constant $C>0$ such that, for any
$\e>0$, any
 solutions $(p,w, q)$ of system (\ref{1202-s1}) satisfies
  \bel{1202-a2}\ba{ll}\ds ||p||_{H^1(Y)}+||w||_{L^2(Y)}+||q||_{L^2(Y)}\\
  \ns\ds\le
 Ce^{C/\e}\big(||p^0||_{L^2(X)}+||w^0||_{L^2(X)}+||q||_{L^2(X^*)}+||q_s||_{L^2(X^*)}\big)\\
 \ns\ds\q+Ce^{-2/\e}\big(||p||_{L^2(X)}+||p_s||_{L^2(X)}+||w||_{L^2(X)}+||q||_{L^2(X)}).
  \ea\ee
 \et

 \ms

{\it Proof of Theorem \ref{1202-t1}. } The proof is based on the global
Carleman estimates presented in Lemma \ref{car1} and Lemma \ref{1202-car1}. In case of $\a=1$, the damping we imposed related to  $d(x)p_s$, in this situation, the main difficulty is to estimate the energy of the coupled system $(p, w, q)$ only localized in $\o_c\cap\o_d$ by $\ds \int_{\o_c\cap\o_d}(|p|^2+|p_s|^2)dx$. The proof is long, hence we divide it
into several steps.

\ss

{\it Step 1. }  Note that there
is no boundary conditions for $p,\  w$ and $q$ at $s=\pm2$ in system
(\ref{1202-s1}). Therefore, we need to introduce a cut-off function
$\varphi=\varphi(s)\in C_0^\infty(-b,b)\subset C_0^\infty(\dbR)$
such that
 \bel{1202-as7}\left\{\ba{ll}\ds
 0\le\varphi(s)\le1 \q &|s|<b,\\
 \ns\ds \varphi(s)=1,&|s|\le b_0,
  \ea\right.\ee
where $1<b_0<b\le2$ are given as follows:
    \bel{1202-as2}
b\=\sqrt{1+\frac{1}{\mu}\ln (2+e^\mu)},\qq
b_0\=\sqrt{b^2-{1\over\mu}\ln \left(\frac{1+e^\mu}{e^\mu}\right)},\qq\forall \mu>\ln 2.
 \ee
Put
 \bel{1202-as8}
\hat p=\varphi p,\q \hat w=\varphi w, \q \hat q=\varphi q.
 \ee
Then, noting that $\varphi$ does not depend on $x$, by (\ref{1202-s1}),
it follows
 \bel{1202-as9}\left\{\ba{ll}\ds
\hat p_{ss}+\D \hat p=F_1& \mbox{ in } X, \\
\ns\ds \hat w_{s}+\D \hat w=F_2&
\mbox{ in } X,\\
\ns\ds \hat q_{s}-\D \hat q=F_3&
\mbox{ in } X,\\
\ns\ds\hat p=\hat q=\hat w=0 &\mbox{ on }\Si,
 \ea\right.
 \ee
where
 \bel{lab-1}\left\{\ba{ll}\ds
 F_1\=\varphi_{ss}p+2\varphi_sp_s+\varphi p^0-i\a d(x)\varphi p_s+c(x)\hat q,\\
 \ns\ds F_2\=\varphi_{s}w+\varphi w^0-i(1-\a)d(x)\varphi q_s+c(x)\hat p,\\
 \ns\ds F_3\=\varphi_{s}q+\hat w.
 \ea\right.\ee

For system (\ref{1202-as9}), by using Lemmas
\ref{car1} and \ref{1202-car1}, noting that
$\a=1$, we conclude that there is a $\mu_1>0$
such that for all $\mu\ge \mu_1$, one can find
two constants $C=C(\mu)>0$ and $\l_1=\l_1(\mu)$
so that  for all $\l\ge\l_1$, it holds that
  \bel{as216}\ba{ll}\ds
 \l\mu^2\int_{-b}^b\int_{\O}\th^2\phi(|\n\hat
 p|^2+|\hat p_s|^2+\l^2\mu^2\phi^2|\hat p|^2+\l^2\mu^2|
 \hat w|^2+\l^2\mu^2|\hat q|^2)dxds\\
 \ns\ds\q+\int_{-b}^b\int_{\O}\th^2(\l\phi)^{-1}(|\D \hat w|^2+|\hat w_s|^2+|\D \hat q|^2+|\hat q_s|^2)dxds\\
\ns\ds\le C\Big[\int_{-b}^b\int_{\O}\th^2\(|F_1|^2+|F_2|^2+|F_3|^2\)dxds\\
\ns\ds\q+C\l\mu^2\int_{-b}^b\int_{\o_1}\th^2\phi \(\l^2\mu^2\phi^2|\hat
p|^2+|\hat p_s|^2+|\n\hat p|^2+\l^2\mu^2\phi^2|\hat w|^2+\l^2\mu^2\phi^2|\hat q|^2\) dxds.
 \ea\ee

 {\it Step 2. } Let us estimate ``$\ds\int_{-b}^b\int_{\o_1}\th^2\phi^3|\hat w|^2dxds$".

 Recall that $\eta_2\in C_0^\i(\o_2) $ satisfying $\eta_2=1$ in $\o_1$. By (\ref{1202-as9}) and (\ref{lab-1}), we have
 \bel{1201-x0}\ba{ll}&\ds
\th^2\phi^3\eta_2^2|\hat w|^2\\
\ns&\ds=\th^2\phi^3\eta_2^2\overline{\hat w}(\hat q_s-\D\hat q)-\th^2\phi^3\eta_2^2\overline{\hat w}\varphi_s q\\
\ns&\ds=-\th^2\phi^3\eta_2^2\hat q\overline{(\hat w_s+\D\hat w)}+(\th^2\phi^3\eta_2^2\overline{\hat w}\hat q)_s-(\th^2\phi^3\eta_2^2)_s\overline{\hat w}\hat q-\th^2\phi^3\eta_2^2\overline{\hat w}\varphi_s q\\
\ns&\ds\q-\sum_{j=1}^n\[\th^2\phi^3\eta_2^2(\overline{\hat w}\hat q_{x_j}-\overline{\hat w}_{x_j}\hat q)\]_{x_j}+\sum_{j=1}^n(\th^2\phi^3\eta_2^2)_{x_j}({\overline {\hat w}}\hat q_{x_j}-{\overline {\hat w}}_{x_j}\hat q).
 \ea\ee
Now, integrating (\ref{1201-x0}) on $(-b,b)\t\O$, noting that  $\hat
w(-b)=\hat w(b)=0$ in $\O$, by (\ref{1202-as9}) and (\ref{1202-as6}), we find that
 \bel{1219-x0}\ba{ll}\ds
 \int_{-b}^b\int_{\O}\th^2\eta_2^2\phi^3|\hat w|^2dxds\\
 \ns\ds\le -\int_{-b}^b\int_{\O}\th^2\phi^3\eta_2^2\hat q {\overline {F_2}}dxds+C\l^2\mu^2\int_{-b}^b\int_{\O}\th^2\phi^5\eta_2^2|\hat q|^2dxds\\
 \ns\ds\q+C\int_{-b}^b\int_{\O}\|\sum_{j=1}^n(\th^2\phi^3\eta_2^2)_{x_j}({\overline {\hat w}}\hat q_{x_j}-{\overline {\hat w}}_{x_j}\hat q)\| dxds.
  \ea\ee
 However, by  Lemma \ref{le0m-01}, it is easy to check that
  \bel{1219-x1}\ba{ll}\ds
 \int_{-b}^b\int_{\O}\|\sum_{j=1}^n(\th^2\phi^3\eta_2^2)_{x_j}({\overline {\hat w}}\hat q_{x_j}-{\overline {\hat w}}_{x_j}\hat q)\| dxds\\
  \ns\ds\le \frac{1}{(\l\mu)^3}\int_{-b}^b\int_{\o_1}\th^2\eta_2^2\phi(|\n \hat w|^2+\phi^2|\hat w|^2)dxds\\
  \ns\ds\q+C(\l\mu)^5\[\int_{-b}^b\int_{\o_1}\th^2\eta_2^2\phi^5|\n \hat q|^2ddxds+\int_{-b}^b\int_{\o_2}\th^2\phi^7|\hat q|^2dxds\]\\
  \ns\ds\le \frac{1}{(\l\mu)^5}\int_{-b}^b\int_{\O}\th^2|F_2|^2dxds+\frac{C}{\l\mu}\int_{-b}^b\int_{\O}\th^2\phi^3\eta_2^2|\hat w|^2dxds\\
  \ns\ds\qq+\frac{C}{(\l\mu)^4}\int_{-b}^b\int_{\O}\th^2|F_3|^2dxds+C(\l\mu)^{14}\int_{-b}^b\int_{\o_2}\th^2\phi^7|\hat q|^2dxds.
  \ea\ee

 Combining (\ref{1219-x0}) and (\ref{1219-x1}), for fixed $\mu$, we conclude that
there is a $\l_3>0$ such that for all $\l\ge \l_3$, one can find
two constants $C=C(\mu)>0$ and $\l_3=\l_3(\mu)$ so that  for all
$\l\ge\l_3$, it holds that
 \bel{1219-x2}\ba{ll}\ds
 \l^3\mu^4\int_{-b}^b\int_{\O}\th^2\eta_2^2\phi^3|\hat w|^2dxds\\
 \ns\ds\le C\[\int_{-b}^b\int_{\O}\th^2\(|F_2|^2+|F_3|^2\)dxds+(\l\mu)^{18}\int_{-b}^b\int_{\O}\th^2\eta_2^2\phi^7|\hat q|^2dxds\].
 \ea\ee

{\it Step 3. } Let us estimate
``$\ds\int_{-b}^b\int_{\o_2}\th^2\phi^7 |\hat q|^2dxds$".

Recall  cut-off function $\eta_3\in C_0^\i(\o_3) $ satisfying $\eta_3=1$ in $\o_2$. By elementary calculation, we have
  \bel{eq2}\ba{ll}\ds
\th^2\phi^7\eta_3^2\overline {\hat q}(\hat
p_{ss}+\D\hat p)\\
 \ns\ds=(\th^2\phi^7\eta_3^2\overline {\hat q}\hat
p_{s})_s-\th^2\phi^7\eta_3^2\hat p_s\overline{\hat
q}_{s}-(\th^2\phi^7\eta_3^2)_s\overline {\hat q}\hat
p_{s}+\th^2\phi^7\eta_3^2\hat p\D\overline{\hat q}\\
\ns\ds\q+\sum_{j=1}^n\[\th^2\phi^7\eta_3^2(\overline {\hat
q}\hat p_{x_j}-\overline {\hat q}_{x_j}\hat p)\]_{x_j}-\sum_{j=1}^n(\th^2\phi^7\eta_3^2)_{x_j}(\overline
{\hat q}\hat p_{x_j}-\overline {\hat q}_{x_j}\hat p).
  \ea \ee

On the other hand, multiplying the first equation of (\ref{1202-as9}) by
$ \th^2\phi^7\eta_3^2\overline {\hat q}$, we have
  \bel{eq3}\ba{ll}\ds c(x)\th^2\phi^7\eta_3^2|\hat q|^2&\ds=\th^2\phi^7\eta_3^2\overline {\hat q}\(\hat
p_{ss}+\D\hat p\)-\th^2\phi^7\eta_3^2\overline {\hat
q}(F_1-c(x)\hat q).
  \ea \ee
Now, integrating (\ref{eq3}) on $(-b,b)\t\O$, noting that  $\hat
p(-b)=\hat p(b)=0$ in $\O$, by (\ref{1202-as9}), (\ref{1202-as6}), (\ref{eq2})
and (\ref{cond1}), we find that
 \bel{eq5}\ba{ll}\ds
\int_{\O}\th^2\eta_3^2\phi^7|\hat q|^2 dxds\\
\ns\ds\le C\int_{-b}^b\int_{\O}\th^2\phi^7|F_1-c(x)\hat q|^2dxds+\frac{C}{(\l\mu)^{16}}\int_{-b}^b\int_{\O}\th^2\phi^{5}\eta_3^2|\n\hat
q|^2dxds\\
\ns\ds\q+\frac{1}{(\l\mu)^{18}}\int_{-b}^b\int_{\O}\th^2(\l\phi)^{-1}\(|\hat q_s|^2+|\D\hat q|^2\)dxds\\
\ns\ds\q+C(\l\mu)^{19}\int_{-b}^b\int_{\o_3}\th^2\phi^{15}(|\n\hat
p|^2+|\hat p_s|^2+|\hat p|^2)dxds.
  \ea\ee
  On the one hand, by Lemma \ref{le0m-01}, we know that
  \bel{122-y}\ba{ll}\ds
   \frac{C}{(\l\mu)^{16}}\int_{-b}^b\int_{\O}\th^2\phi^{5}\eta_3^2|\n\hat
q|^2dxds\\
\ns\ds\le \frac{C}{(\l\mu)^{18}}\int_{-b}^b\int_{\O}\th^2|F_3|^2dxds+\frac{C}{(\l\mu)^{14}}\int_{\O}\th^2\eta_3^2\phi^7|\hat q|^2 dxds.
  \ea\ee
On the other hand, by Lemma \ref{1202-car1},  we have
  \bel{122-y1}\ba{ll}\ds
 \frac{1}{(\l\mu)^{18}}\int_{-b}^b\int_{\O}\th^2(\l\phi)^{-1}\(|\hat q_s|^2+|\D\hat q|^2\)dxds\\
 \ns\ds\le \frac{C}{(\l\mu)^{18}}\[\int_{-b}^b\int_{\O}\th^2|F_3|^2dxds+\l^3\mu^4\int_{-b}^b\int_{\o_1}\th^2\phi^3|\hat q|^2dxds\] \\
  \ns\ds\le \frac{C}{(\l\mu)^{18}}\int_{-b}^b\int_{\O}\th^2|F_3|^2dxds+\frac{C}{(\l\mu)^{14}}\int_{-b}^b\int_{\o_2}\th^2\phi^7|\hat q|^2dxds.
  \ea\ee
 By  (\ref{eq5})--(\ref{122-y1}), for fixed $\mu$, we conclude that
there is a $\l_4>0$ such that for all $\l\ge \l_4$, one can find
two constants $C=C(\mu)>0$ and $\l_4=\l_4(\mu)$ so that  for all
$\l\ge\l_3$, it holds that
 \bel{122-y3}\ba{ll}\ds
 (\l\mu)^{18}\int_{-b}^b\int_{\o_2}\th^2\phi^7 |\hat q|^2dxds\3n&\ds\le C\int_{-b}^b\int_{\O}\th^2\phi^7\((\l\mu)^{18}|F_1-c(x)\hat q|^2+|F_3|^2\)dxds\\
 \ns&\ds\q+C(\l\mu)^{37}\int_{-b}^b\int_{\o_3}\th^2\phi^{15}(|\n\hat
p|^2+|\hat p_s|^2+|\hat p|^2)dxds.
 \ea\ee

Finally, we choose cut-off function $\eta_4\in C_0^\i(\o_c\cap\o_{d}) $ such that $\eta_4=1$ in $\o_3$. Multiplying the first equation of (\ref{1202-as9}) by $\th^2\phi^{15}\eta_4^2\overline{\hat p}$, integrating
it on $(-b,b)\t\O$, using integration by parts,  and noting that
$\hat p(-b)=\hat p(b)=0$ in $\O$, by a simple calculation, we have
  \bel{nr3}\ba{ll}\ds
\int_{-b}^b\int_{\o_3}\th^2\phi^{15}(|\n\hat p|^2+|\hat
p_s|^2)dxds\\
\ns\ds\le C\[e^{C\l}\int_{-b}^b\int_{\o_c\cap\;\o_d}|\hat p|^2 dxds+\frac{1}{(\l\mu)^{37}}\int_{-b}^b\int_{\O}\th^2|F_1|^2dxds\].
 \ea\ee

Combing (\ref{as216}) and (\ref{eq5})--
(\ref{nr3}), by (\ref{cond1}) and (\ref{1202-as7}), noting that
 $\hat p=\varphi p,\ \hat w=\varphi w,\ \hat q=\varphi q$, by (\ref{lab-1}) and noting that $\a=1$, taking $\l_*=\max\{\l_j,\ j=0, 1, 2, 3, 4\}$, for any $\l\ge\l_*$, we have
  \bel{equ216}\ba{ll}\ds
 \l\mu^2\int_{-b}^b\int_{\O}\th^2\phi(|\n
 p|^2+|p_s|^2+\l^2\mu^2\phi^2|p|^2+\l^2\mu^2\phi^2|w|^2+\l^2\mu^2\phi^2|q|^2)dxds\\
\ns\ds\le Ce^{C\l}\Big\{\int_{-b}^b\int_{\O}(|
 p^0|^2+|w^0|^2)dxds+\int_{-b}^b\int_{\o_c\cap\;\o_{d}} (|
p|^2+|p_s|^2)dxds\Big\}\\
 \ns\ds\q+C(\l\mu)^{18}\int_{(-b,-b_0)\bigcup (b_0,
b)}\int_{\O}\th^2\phi^7(|
 p|+|p_s|^2+|
 q|+|w|^2)dxds.
 \ea\ee

 \ms

{\it Step 4. } Recalling (\ref{as4}) and (\ref{1202-as2}) for the definitions of
$\phi$ and $b,\ b_0$, respectively, it is easy to check that
 \bel{as5}\left\{\ba{ll}\ds
 \phi(s,\cd)\ge 2+e^\mu , &\hb{ for any $s$ satisfying }|s|\le 1,\\
 \ns\ds \phi(s,\cd)\le 1+e^{\mu}, &\hb{ for any $s$ satisfying }b_0\le |s|\le
 b.
 \ea\right.\ee
Finally, denote $c_0=2+e^{\mu}>1$. Fixing the parameter $\mu$ in
(\ref{as216}), and using (\ref{as5}), one finds that
 \bel{as217}\ba{ll}\ds
\l e^{2\l c_0}\int_{-1}^1\int_{\O}(|\n p|^2+|p_s|^2+|p|^2+|
w|^2+|q|^2)dxds\\[3mm]
\ns\ds\le
Ce^{C\l}\left\{\int_{-2}^2\int_{\O}(|p^0|^2+|w^0|^2)dxds+\int_{-2}^2\int_{\o_c\cap\;\o_{d}}(|p|^2+|p_s|^2)
dxds\right\}\\[3mm]
\ns\ds\q+C(\l\mu)^{18}e^{c_2\mu} e^{2\l(c_0-1)}\int_{(-b,-b_0)\bigcup (b_0,
b)}\int_{\O}(|p|^2+|p_s|^2+|w|^2+|q|^2)dxds
 \ea\ee
where $c_2=7\mu ||\psi||_{L^{\i}(\O)}>0$.

From (\ref{as217}), one concludes that there exists an $\e_2>0$ such
that the desired inequality (\ref{1202-0a2}) holds for $\e\in (0,\e_2]$,
which, in turn, implies that it holds for any $\e>0$. This completes
the proofs of Theorem \ref{1202-t1}.\endpf

\ms

{\it Proof of Theorem \ref{1202-t2}. }  In case of $\a=0$, the damping we  imposed related to $d(x)q_s$,  in this situation, we will estimate the energy of the coupled system $(p, w, q)$ only localized in $\o_c\cap\o_d$ by $\ds \int_{\o_c\cap\o_d}(|q|^2+|q_s|^2)dx$. We divide the proof into several steps.

\ms

{\it Step 1. } Recall  $\eta_2\in C_0^\i(\o_2)
$ such that $\eta_2=1$ in $\o_1$. Multiplying
the first equation of (\ref{1202-as9}) by
$\th^2\phi\eta_2^2\overline{\hat p}$,
integrating it on $(-b,b)\t\O$, using
integration by parts,  noting that $\hat
p(-b)=\hat p(b)=0$ in $\O$ and $\a=0$, by a
simple calculation, we conclude that
  \bel{0nr3}\ba{ll}\ds
\int_{-b}^b\int_{\O}\eta_2^2\th^2\phi(|\n\hat p|^2+|\hat
p_s|^2)dxds\\
\ns\ds\le C(\l\mu)^2\int_{-b}^b\int_{\O}\th^2\eta_2^2\phi^3|\hat p|^2 dxds+ \frac{1}{(\l\mu)^2}\int_{-b}^b\int_{\O}\th^2|F_1|^2dxds.
 \ea\ee
Therefore, by (\ref{as216}) and (\ref{0nr3}), we conclude that there is a $\mu_1>0$
such that for all $\mu\ge \mu_1$, one can find
two constants $C=C(\mu)>0$ and $\l_1=\l_1(\mu)$
so that  for all $\l\ge\l_1$, it holds that
  \bel{01as216}\ba{ll}\ds
 \l\mu^2\int_{-b}^b\int_{\O}\th^2\phi(|\n\hat
 p|^2+|\hat p_s|^2+\l^2\mu^2\phi^2|\hat p|^2+\l^2\mu^2|
 \hat w|^2+\l^2\mu^2|\hat q|^2)dxds\\
 \ns\ds\le C\int_{-b}^b\int_{\O}\th^2\(|F_1|^2+|F_2|^2+|F_3|^2\)dxds\\
\ns\ds\qq+C\l^3\mu^4\int_{-b}^b\int_{\o_2}\th^2\phi^3(|\hat
p|^2+|\hat w|^2)dxds+Ce^{C\l}\int_{-b}^b\int_{\o_2}|\hat q|^2dxds.
 \ea\ee

\ms

{\it Step 2. } Let us estimate ``$\ds \int_{-b}^b\int_{\o_2}\th^2\phi^3|\hat
p|^2dxds$".

Recall that  $\eta_3\in C_0^\i(\o_3) $ such that
$\zeta_3=1$ in $\o_2$. Then, multiplying the
second equation of (\ref{1202-as9}) by $
\th^2\phi^3\eta_3^2\overline {\hat p}$, we have
\bel{03eq3}\ba{ll}\ds
c(x)\th^2\phi^3\eta_3^2|\hat p|^2
=\th^2\phi^3\eta_3^2\overline {\hat p}\big(\hat
w_{s}+\D\hat w\big)
-\th^2\phi^3\eta_3^2\overline {\hat
p}\big(F_2-c(x)\hat p\big).
  \ea \ee
However,
 \bel{122-f1}\ba{ll}\ds
 \th^2\phi^3\eta_3^2\overline {\hat p}\big(\hat
w_{s}+\D\hat w\big) \\
\ns\ds=(\th^2\phi^3\eta_3^2\overline {\hat
p}\hat w)_{s}-(\th^2\phi^3\eta_3^2)_s\overline
{\hat p}\hat w-\th^2\phi^3\eta_3^2\overline
{\hat p}_s\hat
w\\
\ns\ds\q+\sum_{j=1}^n\big(\th^2\phi^3\eta_3^2\overline
{\hat p}\hat
w_{x_j}\big)_{x_j}-\sum_{j=1}^n(\th^2\phi^3\eta_3^2)_{x_j}\overline
{\hat p}\hat
w_{x_j}-\sum_{j=1}^n\th^2\phi^3\eta_3^2\overline
{\hat p}_{x_j}\hat w_{x_j}.
 \ea\ee

Now, integrating (\ref{03eq3}) on $(-b,b)\t\O$,
by (\ref{122-f1}) and Lemma \ref{1202-car1}, we get that
 \bel{03eq5}\ba{ll}\ds
\int_{-b}^b\int_{\O}\th^2\eta_3^2\phi^3|\hat p|^2 dxds\\
\ns\ds\ds\le C\int_{-b}^b\int_{\O}\th^2\phi^3|F_2-c(x)\hat p|^2 dxds+C(\l\mu)^3\int_{-b}^b\int_{\O}\th^2\eta_3^2\phi^5(|\hat w|^2+|\n \hat w|)^2 dxds\\
\ns\ds\q+\frac{C}{(\l\mu)^3}\int_{-b}^b\int_{\O}\th^2\eta_3^2\phi(|\hat p_s|^2+\n \hat p|^2) dxds.
 \ea\ee
  Proceeding exactly as (\ref{0nr3}), we have
  \bel{122-f2}\ba{ll}\ds
  \frac{C}{(\l\mu)^3}\int_{-b}^b\int_{\O}\th^2\eta_3^2\phi(|\hat p_s|^2+\n p|^2) dxds\\
  \ns\ds\le  \frac{C}{\l\mu}\int_{-b}^b\int_{\O}\th^2\eta_2^2\phi^3|\hat p|^2 dxds+\frac {C}{(\l\mu)^5}\int_{-b}^b\int_{\O}\th^2|F_1|^2dxds.
  \ea\ee
Combining (\ref{03eq5}) and (\ref{122-f2}), for fixed $\mu$, there is a $\l_5>0$ such that for any $\l\ge\l_5$, the following holds:
 \bel{003eq5}\ba{ll}\ds
\int_{-b}^b\int_{\o_2}\th^2\phi^3|\hat p|^2 dxds
&\ds\le C\int_{-b}^b\int_{\O}\th^2\phi^3\big(|F_2-c(x)\hat p|^2+(\l\mu)^{-5}|F_1|^2\big) dxds\\
\ns&\ds\qq+C(\l\mu)^3\int_{-b}^b\int_{\O}\th^2\eta_3^2\phi^5(|\hat
w|^2+|\n \hat w|)^2 dxds.
 \ea\ee

  \ms

 {\it Step 3. } By using Lemma (\ref{le0m-01}) and proceeding similarly analysis as  (\ref{1219-x2}), we have
  \bel{122-fx0}\ba{ll}\ds
  C(\l\mu)^3\int_{-b}^b\int_{\O}\th^2\eta_3^2\phi^5(|\hat w|^2+|\n\hat  w|)^2 dxds\\
  \ns\ds \le C(\l\mu)^6\int_{-b}^b\int_{\O}\th^2\eta_3^2\phi^5|\hat w|^2dxds+C\int_{-b}^b\int_{\O}\th^2\phi^3|F_2|^2 dxds\\
\ns\ds\le  C(\l\mu)^{-4}\int_{-b}^b\int_{\O}\th^2\phi^3\(|F_2|^2+|F_3|^2\)dxds+Ce^{C\l}\int_{-b}^b\int_{\o_c\cap\o_d}|\hat q|^2dxds.
 \ea\ee

Finally, by (\ref{01as216}),  (\ref{1201-x0})
and  (\ref{03eq5}),  noting that $\hat
p=\varphi p,\ \hat w=\varphi w,\ \hat q=\varphi
q$, by (\ref{lab-1}), taking $\l^*=\max\{\l_1,\l_5\}$, for any $\l\ge\l^*$, we obtain that
  \bel{03equ216}\ba{ll}\ds
 \l\mu^2\int_{-b}^b\int_{\O}\th^2\phi(|\n
 p|^2+|p_s|^2+\l^2\mu^2\phi^2|p|^2+\l^2\mu^2\phi^2|w|^2+\l^2\mu^2\phi^2|q|^2)dxds\\
\ns\ds\le Ce^{C\l}\Big[\int_{-b}^b\int_{\O}(|
 p^0|^2+|w^0|^2)dxds+\int_{-b}^b\int_{\o_c\cap\o_{d}}( |
q|^2+|q_s|^2)dxds\Big]\\
 \ns\ds\q+C(\l\mu)^5\int_{(-b,-b_0)\bigcup (b_0,
b)}\int_{\O}\th^2\phi^3(|
 p|+|p_s|^2+|
 q|+|w|^2)dxds.
 \ea\ee
Then, proceeding exactly the same analysis as (\ref{as5}) and (\ref{as217}), we complete the proof of Theorem \ref{1202-t2}.\endpf

\ms

\section{ Proof of the main results}\label{ss6}

In this section, we shall give the proof of
logarithmic decay results. Recall (\ref{0a2})
for the  definition of $\cA$ and $D(\cA)$.
Denote by $\cR(\cA)$ the resolvent set of $\cA$.
The proof of the logarithmic decay relies on a
result of Burq (\cite[Theorem 3]{B}), which
links it to estimates on the resolvent of
$(\cA-i\mu)$ for $\mu\in\dbR$. Later, Duyckaerts
consider the resolvent operator $(\cA-\l
I)^{-1}$ for complex number $\l$. Concerning the
resolvent estimate established in Theorem
\ref{0t2}, let us recall the following known
result.
 \bl\label{Duy}(\cite[Theorem 7.1]{D})
 Let $D>0$, and
  $$
  \cO_D\=\big\{\l\in\dbC: \ |\Re\l|< D^{-1}e^{-D|\Im\l|}\big\}.
  $$
  Assume that for some $D>0$,\ $\cO_D$ is included in $\cR(\cA)$, and that in $\cO_D$ there is a positive constant $C$ such that
   $$
   ||(\cA-\l I)^{-1}||_{\cL(H)}\le Ce^{C|\Im\l|}.
   $$
   Then for all $k$ there exists $C_k$ such that
    $$
    ||e^{t\cA}U_0||_{H}\le \frac{C_k}{(\log (t+2))^k}||U_0||_{D(\cA^k)},  \q \forall U_0\in D(\cA^k).
    $$
  \el

Therefore, once we prove the existence and the
estimate of the norm of the resolvent $(\cA-\l
I)^{-1}$ when $
 \Re\l\in\[-e^{-C|\Im\l|}/ C,0\]$, stated in Theorem \ref{0t2},  by virtue of Lemma \ref{Duy}, we can get Theorem \ref{0t1}, immediately.

\ms

{\it Proof of Theorem \ref{0t2}. } We divide the proof into two
steps.

\ms

{\it Step 1.}  First, fix $F=(f^0,f^1,g^0,g^1)\in H$ and
$U_0\=U(0)=(y^0,y^1,z^0,z^1)\in D(\cA)$. It is easy to see that the
following equation
 \bel{6a1}
(\cA-\l I)U_0=F
 \ee
is equivalent to
 \bel{6a2}\left\{\ba{ll}\ds
 -\l y^0+y^1=f^0, & \hb { in } \O,\\
 \ns\ds
 \D y^0-[\a d(x)+\l] y^1-c(x)z^0=f^1, & \hb { in } \O,\\
 \ns\ds -\l z^0+z^1=g^0, & \hb { in } \O,\\
 \ns\ds
-\D^2 z^0-[(1-\a)d(x)+\l)]z^1-c(x)y^0=g^1,  & \hb { in }
 \O,\\
\ns\ds y^0=z^0=\D z^0=0 &\hb{ on } \G.\\
 \ea\right.\ee
 By (\ref{6a2}), we conclude that
 \bel{6a3}\left\{\ba{ll}\ds
\D y^0 -\l^2 y^0-\a\l (x)y^0-c(x)z^0=[d(x)+\l] f^0+f^1 & \hb { in } \O,\\
\ns\ds-\D^2z^0-\l^2 z^0-(1-\a)\l d(x)z^0-c(x)y^0=\l g^0+g^1 & \hb { in } \O,\\
\ns\ds y^0=z^0=\D z^0=0 &\hb{ on } \G,\\
\ns\ds y^1=f^0+\l y^0,\ z^1=g^0+\l z^0 & \hb { in } \O.
 \ea\right.\ee
Put
 \bel{v}
p=e^{i\l s}y^0,\qq q=e^{i\l s}z^0.
 \ee
It is easy check that $p$ and $q$ satisfy the following equation:
 \bel{6a4}\left\{\ba{ll}\ds
p_{ss}+\D p+i\a d(x) p_s-c(x)q=(\l f^0+df^0+f^1)e^{i\l s} & \hb { in } \dbR\t \O,\\
\ns\ds q_{ss}-\D^2 q+i(1-\a)d(x) q_s-c(x)p=(\l g^0+g^1)e^{i\l s} & \hb { in } \dbR\t \O,\\
\ns\ds p=q=\D q=0& \hb { on }\dbR\t \G.
 \ea\right.\ee

 Further, we set
   \bel{1201-v}
 w=q_s-\D q.
 \ee
 Then, clearly, $p, q$ and $w$ satisfy the following equation:
 \bel{1202-6a4}\left\{\ba{ll}\ds
p_{ss}+\D p+i\a d(x) p_s-c(x)q=(\l f^0+df^0+f^1)e^{i\l s} & \hb { in } \dbR\t \O,\\
\ns\ds w_{s}+\D w+i(1-\a)d(x) q_s-c(x)p=(\l g^0+g^1)e^{i\l s} & \hb { in } \dbR\t \O,\\
\ns\ds q_{s}-\D q=w & \hb { in } \dbR\t \O,\\
\ns\ds p=q=w=0& \hb { on }\dbR\t \G.
 \ea\right.\ee

\ms

{\it Step 2. } By (\ref{v}), we have the
following estimates. \bel{6a7}\left\{\ba{ll}\ds
||y^0||_{H^1_0(\O)}+||z^0||_{H^2(\O)}\le Ce^{C|\Im\l|}\big(||p||_{H^1(Y)}+||w||_{L^2(Y)}+||q||_{L^2(Y)}\big)\\
\ns\ds ||p||_{H^1(X)}+||w||_{L^2(X)}+||q||_{L^2(X)}\le C(|\l|+1)e^{C|\Im\l|}\big(||y^0||_{H^1_0(\O)}+||z^0||_{H^2(\O)}\big),\\
 \ns\ds ||p||_{L^2(X^*)}+||p_s||_{L^2(X^*)}\le C(1+|\l|)e^{C|\Im\l|}||y^0||_{L^2(\o_c\cap\o_d)},\\
 \ns\ds ||q||_{L^2(X^*)}+||q_s||_{L^2(X^*)}\le C(1+|\l|)e^{C|\Im\l|}||z^0||_{L^2(\o_c\cap\o_d)}.
 \ea\right.\ee

 Now, in case of $\a=1$, applying Theorem \ref{1202-t1} to equation (\ref{6a4}), and combining
(\ref{6a7}), we get that
 \bel{6a5}\ba{ll}\ds
||y^0||_{H^1_0(\O)}+||z^0||_{H^2(\O)}\\
\ns\ds\le
Ce^{C|\Im\l|}\(||f^0||_{H^1_0(\O)}+||f^1||_{L^2(\O)}+||g^0||_{H^2(\O)}+||g^1||_{L^2(\O)}+||y^0||_{L^2(\o_c\cap\;\o_d)}\).
 \ea\ee
 In case of $\a=0$, applying Theorem \ref{1202-t2} to equation (\ref{6a4}), and combining
(\ref{6a7}), we obtain that
 \bel{06a5}\ba{ll}\ds
||y^0||_{H^1_0(\O)}+||z^0||_{H^2(\O)}\\
\ns\ds\le
Ce^{C|\Im\l|}\(||f^0||_{H^1_0(\O)}+||f^1||_{L^2(\O)}+||g^0||_{H^2(\O)}+||g^1||_{L^2(\O)}+||z^0||_{L^2(\o_c\cap\;\o_d)}\).
 \ea\ee
On the other hand,
multiplying the first equation of (\ref{6a3}) by $2\overline
y^0$ and integrating it on $\O$, it follows that
 \bel{6a8}\ba{ll}\ds
\int_{\O}\(-\D y^0+\l^2y^0+\a\l d(x)y^0+c(x)z^0\)\cd2\overline y^0dx\\
\ns\ds=2\l^2\int_{\O}|y^0|^2dx+2\int_\O |\n y^0|^2dx+2\a\int_{\O}\l d(x) |y^0|^2dx+2\int_\O c(x)
z^0\overline y^0dx.
 \ea\ee
Similarly, multiplying the second equation of
(\ref{6a3}) by $2\overline z^0$ and integrating
it on $\O$, we find that \bel{n6a8}\ba{ll}\ds
\int_{\O}\(\D^2 z^0+\l^2z^0+(1-\a)\l
d(x)z^0+cy^0\)\cd2\overline
z^0dx\\
\ns\ds=2\l^2\int_\O|z^0|^2dx+2\int_\O
|\D z^0|^2dx+2(1-\a)\int_{\O}\l d(x) |z^0|^2dx+2\int_\O c(x)
y^0\overline z^0dx.
 \ea\ee

Taking the imaginary part in the both sides of (\ref{6a8})
+(\ref{n6a8}), by (\ref{6a4}), we have
 \bel{6a9}\ba{ll}\ds
2\a|\Im\l|\int_{\O}d(x)|y^0|^2dx+2(1-\a)|\Im\l|\int_{\O}d(x)|z^0|^2dx\\
\ns\ds \le C\[||(\l
f^0+df^0+f^1)||_{L^2(\O)}||y^0||_{L^2(\O)}+|\Im\l||\Re\l|||y^0||^2_{L^2(\O)}\]\\
\ns\ds\q+C\[||(\l
g^0+g^1)||_{L^2(\O)}||z^0||_{L^2(\O)}+|\Im\l||\Re\l||| z^0||^2_{L^2(\O)}\],
 \ea\ee
where we have use the following obvious fact:
 $$
\Im\int_\O c(x) z^0\overline y^0dx+\Im\int_\O c(x)
y^0\overline z^0dx=\Im \int_\O c(x) (z^0\overline
y^0+y^0\overline z^0)dx=0.
 $$

Hence, combining (\ref{6a5}) and (\ref{6a9}), and noting that
$d(x)\ge d_0>0$ on $\o_{d}$, we arrive at
 \bel{6a10}\ba{ll}\ds
||y^0||_{H^1_0(\O)}+||z^0||_{H^2(\O)}\\
\ns\ds\le
Ce^{C|\Im\l|}\big(||f^0||_{H^1_0(\O)}+||f^1||_{L^2(\O)}+||g^0||_{H^2(\O)}+||g^1||_{L^2(\O)}\big)\\
\ns\ds\q+Ce^{C|\Im\l|}|\Re\l|(||y^0||_{H^1_0(\O)}+||z^0||_{H^2(\O)}).
 \ea\ee
We now take
 $$
Ce^{C|\Im\l|}|\Re\l|\le\frac{1}{2},
$$
which holds, whenever $\ds |\Re\l|\le e^{-C|\Im\l|}/C$ for some
sufficiently large $C>0$. Then, by (\ref{6a10}), we find that
  \bel{6a11}\ba{ll}\ds
||y^0||_{H^1_0(\O)}+||z^0||_{H^2(\O)}\\
\ns\ds\le
Ce^{C|\Im\l|}\big(||f^0||_{H^1_0(\O)}+||f^1||_{L^2(\O)}+||g^0||_{H^2(\O)}+||g^1||_{L^2(\O)}\big).
 \ea\ee
Recalling that $y^1=f^0+\l y^0,\ z^1=g^0+\l z^0$, it follows
  \bel{6a12}\ba{ll}\ds
||y^1||_{L^2(\O)}+||z^1||_{L^2(\O)}\\
\ns\ds\le
||f^0||_{L^2(\O)}+|\l|||y^0||_{L^2(\O)}+||g^0||_{H^2(\O)}+|\l|||z^0||_{H^2(\O)}\\
\ns\ds\le
Ce^{C|\Im\l|}\big(||f^0||_{H^1_0(\O)}+||f^1||_{L^2(\O)}+||g^0||_{H^2(\O)}+||g^1||_{L^2(\O)}\big).
 \ea\ee
By (\ref{6a11})--(\ref{6a12}), we know that $\cA-\l I$ is injective.
Therefore $\cA-\l I$ is bi-injective from $D(\cA)$ to $H_1$.
Moreover,
 $$
||(\cA-\l I)^{-1}||_{\cL(H;H)}\le
Ce^{C|\Im\l|},\q  \Re\l\in
[-e^{-C|\Im\l|}/C,0],\qq |\l|\ge 1.
 $$
This completes the proof of Theorem \ref{0t2}.\endpf

\ms

\end{document}